\documentclass[11pt]{article}
\usepackage{amsmath,amssymb,amsthm,mathrsfs}
\usepackage[margin=1.5in]{geometry}
\usepackage[hidelinks]{hyperref}

\newtheorem{theorem}{Theorem}
\newtheorem{proposition}{Proposition}
\newtheorem{lemma}{Lemma}
\newtheorem{corollary}{Corollary}
\theoremstyle{definition}
\newtheorem{definition}{Definition}
\newtheorem{remark}{Remark}

\newcommand{\R}{\mathbb{R}}
\newcommand{\N}{\mathbb{N}}
\newcommand{\dd}{\,\mathrm{d}}
\newcommand{\Poly}{\mathcal{P}}
\newcommand{\sob}{W^{\alpha,2}_\rho(\R)}
\DeclareMathOperator{\sech}{sech}
\DeclareMathOperator{\sign}{sign}

\begin{document}

\title{Suboptimality of Gaussian quadrature for Freud and logistic weights in weighted Sobolev spaces}

\author{Yuya Suzuki\\[1ex]
  {\small yuya.suzuki@bageri.fi, Aalto University, Finland}\\
}

\maketitle

\begin{abstract}
We study the worst-case error of Gaussian quadrature for integrals over the real line weighted by Freud weights $\exp(-|x|^\beta)$ with $\beta>1$, and by the logistic density $\mathrm{e}^{-x}/(1+\mathrm{e}^{-x})^2$. The integrand functions are assumed to live in the $\rho$-weighted $L^2$-based Sobolev space of smoothness $\alpha\in\N$.

 We prove that any linear quadrature using the $n$ Gaussian nodes, irrespective of the
choice of quadrature weights, has a worst-case error bounded from below by a constant
times $n^{-\alpha(1-1/\beta)}$ for the Freud weights, and by a constant times
$(\log n)^{-\alpha}$ for the logistic weight. Our proof is based on three simple
observations. First, the worst-case error of any linear quadrature is bounded from
below in terms of the number of its nodes in a fixed interval. Second, the number of
Gaussian nodes in such an interval is bounded from above by means of Christoffel
functions. Third, Christoffel functions admit sharp lower bounds on such an interval.
For the Freud weights with even integer exponents and for the logistic weight, all
steps are proved in an elementary and self-contained way; for the logistic weight, we
use explicit recurrence coefficients of the corresponding orthogonal polynomials. Our
results provide a theoretical explanation, for integrands of finite smoothness, of the
slow convergence of Gauss--Logistic quadrature observed in the literature. In contrast,
for the logistic weight and for the Freud weights with even integer exponents, the
M\"obius-transformed trapezoidal rule attains the optimal rate $n^{-\alpha}$. This rule
does not depend on the smoothness $\alpha$ and requires only point values of the
integrand and the weight.
\end{abstract}

\section{Introduction}\label{sec:intro}

This paper considers numerical integration over the real line,
\[
  I_\rho(f):=\int_\R f(x)\,\rho(x)\dd x,
\]
where $\rho$ is a positive weight function and $f$ is finitely smooth. 
A classical and widely used approach is Gaussian quadrature, which uses the zeros of the orthogonal polynomials with respect to $\rho$ as nodes and integrates all polynomials of degree at most $2n-1$ exactly with $n$ nodes; see, e.g., \cite{Gau04,Sze}.
On bounded intervals, Gaussian quadrature is known to be essentially optimal for analytic integrands \cite{Bak67,Pet98}. 
On the real line, however, the situation is different, even for analytic integrands.
Trefethen \cite{Tre22} pointed out that Gauss--Hermite quadrature is less efficient than, e.g., the trapezoidal rule.
This was made rigorous in \cite{GKT24}: in weighted Hardy spaces of analytic functions, the worst-case error of Gauss--Hermite quadrature cannot decay faster than $\exp(-C\sqrt n)$, which is slower than the rates attained by suitably truncated trapezoidal rules and scaled Gauss--Legendre and Clenshaw--Curtis rules.

For integrands with finite smoothness, the suboptimality of Gauss--Hermite quadrature was first shown for Lipschitz continuous functions in \cite{Cur98}. 
When $\rho$ is the standard Gaussian density, the corresponding weighted Sobolev spaces have been considered by numerous authors; see, e.g., \cite{DIL18,DN24,KSG23} and the references therein.
In \cite{KSG23}, it was proved that Gauss--Hermite quadrature attains merely the rate $n^{-\alpha/2}$ in the Gaussian Sobolev space of smoothness $\alpha$, while the optimal rate among all linear quadratures is $n^{-\alpha}$. 
In both \cite{GKT24} and \cite{KSG23}, the lower bounds exploit the distribution of the Gauss--Hermite nodes and are independent of the quadrature weights. 
However, all of these results are restricted to the Gaussian weight, and it has remained unclear how the suboptimality of Gaussian quadrature depends on the weight.

More recently, the M\"obius-transformed trapezoidal rule was introduced in \cite{SHK}.
It was proved to attain the optimal rate $n^{-\alpha}$ for weighted Sobolev spaces with any monotonic Schwartz weight, a class of weights that includes weights decaying only at the rate $\mathrm{e}^{-|x|}$. 
The optimal rate has recently been established also for polynomial weights, in related polynomially weighted Sobolev spaces \cite{HS2026}.
In particular, this applies to the logistic weight and to the Freud weights with even
integer exponents (see Section~\ref{sec:mtr}). For general Freud weights, the optimal
rate $n^{-\alpha}$ follows from a recent result of D\~ung \cite[Corollary~3.3]{Dun25}.
The numerical experiments in \cite{SHK} further showed that Gauss--Logistic quadrature, i.e., Gaussian quadrature for the logistic density, converges much more slowly than the M\"obius-transformed trapezoidal rule, and that this slow convergence is not caused by the limited precision of floating point arithmetic.
Similarly, Mastroianni and Milovanovi\'c observed that Gauss--Logistic quadrature converges
considerably more slowly than Gaussian quadrature on the half-line, even though the
integrands are analytic in a strip around the real line \cite[Table~5.4.8]{MM08}. They
attributed this, tentatively, to the location of the poles of the integrands. 
To the best of our knowledge, no theoretical explanation of this phenomenon has been available.

In this paper, we show that the suboptimality of Gaussian quadrature is not specific to the Gaussian weight, and that it becomes more severe the more slowly the weight decays. We prove lower bounds on the worst-case error of Gaussian quadrature for two classes of weights: the Freud weights $\rho_\beta(x)=\exp(-|x|^\beta)$ with $\beta>1$, and the logistic density $\rho_{\mathrm L}(x)=\mathrm{e}^{-x}/(1+\mathrm{e}^{-x})^2$. More precisely, we show in Theorem~\ref{thm:main} that any linear quadrature using the $n$ Gaussian nodes, with arbitrary quadrature weights, has worst-case error bounded below by
\[
  c\,n^{-\alpha(1-1/\beta)}\quad\text{for }\rho_\beta,
  \qquad\text{and}\qquad
  c\,(\log n)^{-\alpha}\quad\text{for }\rho_{\mathrm L},
\]
with a constant $c>0$ independent of $n$. For $\beta=2$, the first bound recovers the rate $n^{-\alpha/2}$ of \cite{KSG23} for Gauss--Hermite quadrature. The second bound shows that, to attain a worst-case error $\varepsilon$, the number of Gauss--Logistic nodes must grow at least like $\exp(c'\varepsilon^{-1/\alpha})$.

Our argument is simple and consists of three steps. First, the worst-case error of any linear quadrature with $m$ nodes in the interval $(0,1)$ is bounded below by a constant times $(m+1)^{-\alpha}$; this holds for any weight that is positive and continuous on $[0,1]$ (Lemma~\ref{lem:fooling}). 
Second, the number of Gaussian nodes in $[0,1]$ is bounded above by $\rho(\R)/\min_{[0,1]}\lambda_n$, where $\lambda_n$ denotes the Christoffel function (Lemma~\ref{lem:count}). 
Third, the Christoffel function is bounded below on $[0,1]$ by a constant times $n^{1/\beta-1}$ for $\rho_\beta$, and by a constant times $(\log n)^{-1}$ for $\rho_{\mathrm L}$ (Sections~\ref{sec:interval}--\ref{sec:logistic}). 
For the Freud weights with even integer exponents and for the logistic weight, the third step is proved in a self-contained way. 
For general Freud weights it follows from two standard inequalities from the theory of weighted polynomial approximation \cite{LL90,LL01,MS85}. Christoffel functions are a classical tool in the theory of orthogonal polynomials; see, e.g., \cite{Fre71,Nev86} and, for Freud weights, \cite{LL92,LL01,Lub07}. For the logistic weight, we make use of the fact that the corresponding orthogonal polynomials are, up to scaling, continuous Hahn polynomials \cite{KLS10}, so that the three-term recurrence coefficients are explicitly available.

The rest of the present paper is organized as follows. In Section~\ref{sec:prelim}, we introduce the weights, the weighted Sobolev spaces, and the necessary facts on Christoffel functions and Gaussian quadrature. Section~\ref{sec:main} presents our main lower bound, Theorem~\ref{thm:main}, together with its proof. In Section~\ref{sec:mtr}, we show that the M\"obius-transformed trapezoidal rule, with and without modification, attains the optimal rate for all the weights considered. Finally, Section~\ref{sec:concluding} presents concluding remarks.

Throughout this paper, $\N:=\{1,2,\dots\}$ denotes the set of natural numbers and $\R$ the set of real numbers. For $k\in\N\cup\{0\}$, we denote by $\Poly_k$ the space of real polynomials of degree at most $k$. For a Borel set $A\subset\R$, $\mathbf 1_A$ denotes its indicator function. For an open interval $J\subset\R$, $W^{\alpha,2}(J)$ denotes the standard (unweighted) Sobolev space on $J$.

\section{Preliminaries}\label{sec:prelim}

\subsection{Weights and weighted Sobolev spaces}\label{sec:weights}

In this paper, a \emph{weight} is a positive continuous function $\rho\colon\R\to(0,\infty)$ such that $\int_\R|x|^k\rho(x)\dd x<\infty$ for all $k\in\N$. We mainly consider two classes of weights. The first one consists of the \emph{Freud weights}
\[
  \rho_\beta(x):=\exp(-|x|^\beta),\qquad \beta>1,
\]
for which $\int_\R\rho_\beta(x)\dd x=2\,\Gamma(1+1/\beta)$. The second one is the \emph{logistic weight}
\[
  \rho_{\mathrm L}(x):=\frac{\mathrm{e}^{-x}}{(1+\mathrm{e}^{-x})^2}=\frac{1}{4\cosh^2(x/2)},
\]
which is the density of the standard logistic distribution and hence satisfies $\int_\R\rho_{\mathrm L}(x)\dd x=1$. Both $\rho_\beta$ and $\rho_{\mathrm L}$ are even and decreasing on $[0,\infty)$. Moreover, since $4\cosh^2(x/2)=\mathrm{e}^{x}+2+\mathrm{e}^{-x}\ge \mathrm{e}^{|x|}$, we have $\rho_{\mathrm L}(x)\le\mathrm{e}^{-|x|}$. Thus the logistic weight decays exactly at the rate of the borderline Freud weight $\exp(-|x|)$.

We denote by $L^2_\rho(\R)$ the space of Lebesgue measurable functions $f\colon\R\to\R$ with $\int_\R|f(x)|^2\rho(x)\dd x<\infty$. For $\alpha\in\N$, our target functions live in the weighted Sobolev space
\[
  W^{\alpha,2}_\rho(\R):=\Bigl\{f\in L^2_\rho(\R)\;\Big|\;\|f\|_{W^{\alpha,2}_\rho(\R)}:=\Bigl(\sum_{\tau=0}^{\alpha}\int_\R|f^{(\tau)}(x)|^2\rho(x)\dd x\Bigr)^{1/2}<\infty\Bigr\},
\]
where $f^{(\tau)}$ denotes the $\tau$th weak derivative of $f$. Since $\alpha\ge1$ and $\rho$ is positive and continuous, we have $\sob\subset W^{1,2}_{\mathrm{loc}}(\R)$. Hence every $f\in\sob$ has a continuous representative, which we always use, so that point evaluations are well defined.

For a linear quadrature
\[
  A_n(f):=\sum_{j=1}^{n}v_jf(y_j),\qquad v_j\in\R,\ y_j\in\R,
\]
we measure its quality by the worst-case error
\[
  \|I_\rho-A_n\|_{\mathscr L(\sob,\R)}:=\sup_{0\ne f\in\sob}\frac{|I_\rho(f)-A_n(f)|}{\|f\|_{\sob}},
\]
which may be infinite. Note that $I_\rho$ is a bounded linear functional on $\sob$, since $|I_\rho(f)|\le\rho(\R)^{1/2}\|f\|_{L^2_\rho(\R)}$ by the Cauchy--Schwarz inequality, where $\rho(\R):=\int_\R\rho(x)\dd x$.

\begin{remark}[Choice of the norm]\label{rem:norm}
The $L^2_\rho$-based norm is essential for our setting. If one instead uses norms of the form $\sum_\tau\|f^{(\tau)}\rho\|_{L^1(\R)}$ or $\max_\tau\|f^{(\tau)}\rho\|_{L^\infty(\R)}$, then, for any $n$-point linear quadrature, one can take a bump function supported to the right of the largest node such that $f\rho$ is bounded from above and below by fixed positive constants on a fixed interval. Such a function is annihilated by the quadrature, while its integral is bounded away from zero, so that the worst-case error does not decay at all. In contrast, in $\sob$ the integral of a normalized bump function supported in $[R,R+1]$ is at most $\max_{[R,R+1]}\rho^{1/2}$.
\end{remark}

Following \cite{SHK}, we also recall the class of monotonic Schwartz weights. Let $\mathcal S$ be the Schwartz space on $\R$. Then
\[
  \mathcal S^{\mathrm{mon}}_+:=\bigl\{\rho\in\mathcal S\;\big|\;\rho(x)>0\text{ for all }x\in\R,\ \exists K>0\text{ such that }x\,\rho'(x)\le0\text{ for all }x\in\R\setminus[-K,K]\bigr\}.
\]

\subsection{Christoffel functions and Gaussian quadrature}\label{sec:christoffel}

In this subsection, $\nu$ denotes a positive Borel measure on $\R$ with infinite support and finite moments of all orders.

\begin{definition}[Christoffel function]
For $n\in\N$ and $x\in\R$, the $n$th Christoffel function of $\nu$ is defined by
\[
  \lambda_n(\nu;x):=\inf_{q}\Bigl\{\int_\R q(t)^2\dd\nu(t)\;\Big|\;q\in\Poly_{n-1},\ q(x)=1\Bigr\}.
\]
When $\nu(\mathrm{d}t)=\rho(t)\dd t$ for a weight $\rho$, we write $\lambda_n(\rho;x)$.
\end{definition}

Let $(p_k)_{k\ge0}$ be the orthonormal polynomials with respect to $\nu$, with $\deg p_k=k$ and positive leading coefficients. These are obtained by the Gram--Schmidt process applied to $1,t,t^2,\dots$, which is possible since $\nu$ has infinite support. We denote by
\[
  K_n(x,y):=\sum_{k=0}^{n-1}p_k(x)p_k(y)
\]
the Christoffel--Darboux kernel.

The following lemma collects elementary properties of Christoffel functions. They are known in the literature (see, e.g., \cite{Fre71,Nev86}), and we include the proofs for completeness.

\begin{lemma}\label{lem:christoffel}
Let $\nu,\nu'$ be measures as above, $n\in\N$ and $x\in\R$.
\begin{enumerate}
\item[(i)] It holds that $\lambda_n(\nu;x)=1/K_n(x,x)$.
\item[(ii)] It holds that $\lambda_{n+1}(\nu;x)\le\lambda_n(\nu;x)$.
\item[(iii)] If $\nu'(A)\le\nu(A)$ for all Borel sets $A\subset\R$, then $\lambda_n(\nu';x)\le\lambda_n(\nu;x)$.
\item[(iv)] Let $k\in\N\cup\{0\}$ and $G\in\Poly_k$. Suppose that
\[
  \int_\R g(t)\,G(t)^2\dd\nu'(t)\le\int_\R g(t)\dd\nu(t)
\]
for all Borel functions $g\ge0$. Then
\[
  \lambda_n(\nu;x)\ge G(x)^2\,\lambda_{n+k}(\nu';x).
\]
\end{enumerate}
\end{lemma}

\begin{proof}
(i) Any $q\in\Poly_{n-1}$ can be uniquely written as $q=\sum_{k=0}^{n-1}c_kp_k$
with $c_k\in\R$, and by orthonormality
$\int_\R q(t)^2\dd\nu(t)=\sum_{k=0}^{n-1}c_k^2$. Suppose $q(x)=1$. Then the
Cauchy--Schwarz inequality gives
\[
  1=\Bigl(\sum_{k=0}^{n-1}c_kp_k(x)\Bigr)^2
  \le\Bigl(\sum_{k=0}^{n-1}c_k^2\Bigr)\Bigl(\sum_{k=0}^{n-1}p_k(x)^2\Bigr)
  =K_n(x,x)\int_\R q(t)^2\dd\nu(t).
\]
Since $p_0$ is a nonzero constant, we have $K_n(x,x)>0$. Taking the infimum over
all $q\in\Poly_{n-1}$ with $q(x)=1$ gives $\lambda_n(\nu;x)\ge1/K_n(x,x)$.

For the reverse inequality, consider the polynomial $q^*\in\Poly_{n-1}$ defined
by
\[
  q^*(t):=\frac{K_n(x,t)}{K_n(x,x)}
  =\sum_{k=0}^{n-1}\frac{p_k(x)}{K_n(x,x)}\,p_k(t),\qquad t\in\R.
\]
Then $q^*(x)=1$, and by orthonormality
\[
  \int_\R q^*(t)^2\dd\nu(t)
  =\sum_{k=0}^{n-1}\frac{p_k(x)^2}{K_n(x,x)^2}
  =\frac{1}{K_n(x,x)}.
\]
Hence $\lambda_n(\nu;x)\le1/K_n(x,x)$, which completes the proof of (i).

(ii) The infimum defining $\lambda_{n+1}(\nu;x)$ is taken over all $q\in\Poly_n$
with $q(x)=1$. Since $\Poly_{n-1}\subset\Poly_n$, this set contains all
$q\in\Poly_{n-1}$ with $q(x)=1$, over which the infimum defining
$\lambda_n(\nu;x)$ is taken. Hence $\lambda_{n+1}(\nu;x)\le\lambda_n(\nu;x)$.

(iii) For any $q\in\Poly_{n-1}$ we have
$\int_\R q(t)^2\dd\nu'(t)\le\int_\R q(t)^2\dd\nu(t)$, since $q^2\ge0$. The claim
follows by taking the infimum on both sides over all $q\in\Poly_{n-1}$ with
$q(x)=1$.

(iv) Let $q\in\Poly_{n-1}$ with $q(x)=1$. Applying the assumption with $g=q^2$
yields
\[
  \int_\R q(t)^2\dd\nu(t)\ge\int_\R q(t)^2G(t)^2\dd\nu'(t).
\]
If $G(x)=0$, there is nothing to prove, since the right-hand side of the claim
vanishes. Otherwise, define the polynomial
\[
  \tilde q(t):=\frac{q(t)\,G(t)}{G(x)},\qquad t\in\R.
\]
Here $x$ is fixed, so that $G(x)$ is a nonzero constant. Since $\deg q\le n-1$ and
$\deg G\le k$, we have $\tilde q\in\Poly_{n+k-1}$. Moreover, $\tilde q(x)=q(x)=1$.
Hence, by the definition of $\lambda_{n+k}(\nu';x)$,
\[
  \int_\R q(t)^2G(t)^2\dd\nu'(t)=G(x)^2\int_\R\tilde q(t)^2\dd\nu'(t)
  \ge G(x)^2\lambda_{n+k}(\nu';x).
\]
Taking the infimum over all $q\in\Poly_{n-1}$ with $q(x)=1$ completes the proof.
\end{proof}

It is well known that the zeros $x_{n,1}<\dots<x_{n,n}$ of $p_n$ are real and simple, and that there exist unique $w_{n,1},\dots,w_{n,n}\in\R$ such that
\begin{equation}\label{eq:gauss}
  \int_\R q(t)\dd\nu(t)=\sum_{j=1}^n w_{n,j}\,q(x_{n,j})\qquad\text{for all }q\in\Poly_{2n-1};
\end{equation}
see, e.g., \cite{Gau04,Sze}.  We call $x_{n,1},\dots,x_{n,n}$ the \emph{Gaussian nodes} and $w_{n,1},\dots,w_{n,n}$ the \emph{Gaussian quadrature weights}. The corresponding Gaussian quadrature for a weight $\rho$ is denoted by
\[
  Q^{\mathrm G}_{\rho,n}(f):=\sum_{j=1}^nw_{n,j}f(x_{n,j}).
\]
We refer to $Q^{\mathrm G}_{\rho_\beta,n}$ as Gauss--Freud quadrature and to $Q^{\mathrm G}_{\rho_{\mathrm L},n}$ as Gauss--Logistic quadrature.

\begin{lemma}\label{lem:gauss-weights}
For $j=1,\dots,n$, it holds that $w_{n,j}=\lambda_n(\nu;x_{n,j})>0$. Moreover, $\sum_{j=1}^nw_{n,j}=\nu(\R)$.
\end{lemma}

\begin{proof}
Let $\ell_k\in\Poly_{n-1}$ be the Lagrange basis polynomial with $\ell_k(x_{n,i})=\delta_{ki}$. Since $\ell_k^2\in\Poly_{2n-2}$, the exactness \eqref{eq:gauss} gives
\[
  w_{n,k}=\sum_iw_{n,i}\ell_k(x_{n,i})^2=\int_\R\ell_k(t)^2\dd\nu(t)>0 ,
\]
where positivity follows from $\ell_k\not\equiv0$ and the infinite support of $\nu$.

   Now let $q\in\Poly_{n-1}$ with $q(x_{n,j})=1$. Since $q^2\in\Poly_{2n-2}$ and all Gaussian quadrature weights are positive,
\[
  \int_\R q(t)^2\dd\nu(t)=\sum_{i=1}^nw_{n,i}q(x_{n,i})^2\ge w_{n,j}.
\]
Taking the infimum over all $q\in\Poly_{n-1}$ with $q(x_{n,j})=1$ gives $\lambda_n(\nu;x_{n,j})\ge w_{n,j}$. Conversely, $\ell_j\in\Poly_{n-1}$ satisfies $\ell_j(x_{n,j})=1$ and $\int_\R\ell_j(t)^2\dd\nu(t)=w_{n,j}$ by the first part of the proof, so that $\lambda_n(\nu;x_{n,j})\le w_{n,j}$. The last claim follows from \eqref{eq:gauss} with $q\equiv1$.
\end{proof}

\section{Lower bounds for Gaussian quadrature}\label{sec:main}

Our main theorem is as follows.

\begin{theorem}\label{thm:main}
Let $\alpha\in\N$.
\begin{enumerate}
\item[(i)] Let $\beta>1$ and $\rho=\rho_\beta$. Then there exists a constant $c_{\alpha,\beta}>0$, depending only on $\alpha$ and $\beta$, such that, for any $n\in\N$ and any linear quadrature $A_n(f)=\sum_{j=1}^nv_jf(x_{n,j})$ using the Gaussian nodes for $\rho_\beta$ with arbitrary $v_j\in\R$,
\[
  \|I_\rho-A_n\|_{\mathscr L(\sob,\R)}\ge c_{\alpha,\beta}\,n^{-\alpha(1-1/\beta)} .
\]
\item[(ii)] Let $\rho=\rho_{\mathrm L}$. Then there exists a constant $c_\alpha>0$, depending only on $\alpha$, such that, for any $n\ge2$ and any linear quadrature $A_n(f)=\sum_{j=1}^nv_jf(x_{n,j})$ using the Gaussian nodes for $\rho_{\mathrm L}$ with arbitrary $v_j\in\R$,
\[
  \|I_\rho-A_n\|_{\mathscr L(\sob,\R)}\ge c_\alpha\,(\log n)^{-\alpha} .
\]
\end{enumerate}
In particular, the above lower bounds hold for Gauss--Freud quadrature $Q^{\mathrm G}_{\rho_\beta,n}$ and Gauss--Logistic quadrature $Q^{\mathrm G}_{\rho_{\mathrm L},n}$.
\end{theorem}

The proof is given in Section~\ref{sec:proof-main}, after the following preparatory steps. In Section~\ref{sec:fooling}, we bound the worst-case error of an arbitrary linear quadrature from below in terms of the number of its nodes in $(0,1)$. In Section~\ref{sec:count}, we bound this number for Gaussian nodes by the Christoffel function. In Sections~\ref{sec:interval}--\ref{sec:logistic}, we derive lower bounds for the Christoffel functions of $\rho_\beta$ and $\rho_{\mathrm L}$ on $[0,1]$.

\subsection{A localized lower bound for arbitrary linear quadratures}\label{sec:fooling}

For the unweighted integral over an interval, it is well known that any linear
quadrature with $n$ nodes has worst-case error bounded below by a constant times
$n^{-\alpha}$ in the Sobolev space of smoothness $\alpha$, even when restricted
to functions supported in the interval \cite{Tem92}; see also \cite{DN24,SHK}.
Here, we need a localized version of this bound for the weighted integral
$I_\rho$, in which only the nodes in the interval $(0,1)$ are counted. This
localization is essential, since most of the Gaussian nodes lie outside $(0,1)$.
We prove this bound for completeness, using a bump function argument; see, e.g.,
\cite[Section~2.7]{DHP15}.

The idea of the following lemma is to construct a function in $\sob$ which is
supported in $(0,1)$ and vanishes at all nodes, yet has a large integral relative
to its norm. Fix once and for all $\varphi\in C^\infty_c((0,1))$ with $\varphi\ge0$
and $\varphi\not\equiv0$. For instance, one can take the standard mollifier
rescaled to $(\frac14,\frac34)$, that is, $\varphi(t):=\exp\bigl(-1/(1-(4t-2)^2)\bigr)$
for $t\in(\frac14,\frac34)$ and $\varphi(t):=0$ otherwise. Define
\[
  \Phi_0:=\int_0^1\varphi(t)\dd t,\qquad
  \Phi_\alpha:=\sum_{\tau=0}^{\alpha}\int_0^1|\varphi^{(\tau)}(t)|^2\dd t .
\]

\begin{lemma}\label{lem:fooling}
Let $\alpha\in\N$ and let $\rho$ be a weight. Set
$\rho_{\max}:=\max_{x\in[0,1]}\rho(x)$ and $\rho_{\min}:=\min_{x\in[0,1]}\rho(x)$.
Let $A_n(f)=\sum_{j=1}^nv_jf(y_j)$ be a linear quadrature, and let $m$ be the
number of its distinct nodes in $(0,1)$. Then
\[
  \|I_\rho-A_n\|_{\mathscr L(\sob,\R)}\ge C_{\rho,\alpha}\,(m+1)^{-\alpha},
  \qquad
  C_{\rho,\alpha}:=\frac{\rho_{\min}\,\Phi_0}{2^{\alpha+1}(\rho_{\max}\,\Phi_\alpha)^{1/2}},
\]
where $C_{\rho,\alpha}>0$ is independent of $n$, of the nodes $y_j$ and of the
weights $v_j$.
\end{lemma}

\begin{proof}
Set $M:=2(m+1)$ and $h:=1/M$, and let $I_i:=((i-1)h,ih)$, $i=1,\dots,M$. Let
\[
  S:=\{i\in\{1,\dots,M\}\mid I_i\text{ contains no node of }A_n\}.
\]
Since each of the $m$ nodes in $(0,1)$ lies in at most one $I_i$, we have
$\#S\ge M-m\ge M/2$. Define
\[
  f(x):=\sum_{i\in S}\varphi\Bigl(\frac{x-(i-1)h}{h}\Bigr).
\]
The $i$th term is supported in a compact subset of $I_i$, and these supports are
pairwise disjoint. Hence $f\in C^\infty_c((0,1))\subset\sob$ and $f\ge0$.
Moreover, $f$ vanishes outside $\bigcup_{i\in S}I_i$, which contains no node by
the definition of $S$. Hence $A_n(f)=0$ and $|I_\rho(f)-A_n(f)|=I_\rho(f)$.

Since $f\ge0$, $\rho\ge\rho_{\min}$ on $[0,1]$ and $h\,\#S\ge\frac12$, the
substitution $x=(i-1)h+ht$ yields
\[
  I_\rho(f)\ge\rho_{\min}\,\#S\,h\,\Phi_0\ge\frac12\,\rho_{\min}\Phi_0 .
\]
On the other hand, by the chain rule and the disjointness of the supports, for
$\tau=0,\dots,\alpha$,
\[
  \int_\R|f^{(\tau)}(x)|^2\rho(x)\dd x
  \le\rho_{\max}\,\#S\,h^{1-2\tau}\int_0^1|\varphi^{(\tau)}(t)|^2\dd t
  \le\rho_{\max}\,h^{-2\alpha}\int_0^1|\varphi^{(\tau)}(t)|^2\dd t,
\]
where we used $h\,\#S\le1$ and $h\le1$. Summing over $\tau$ gives
$\|f\|_{\sob}^2\le\rho_{\max}\Phi_\alpha\,h^{-2\alpha}$. Since
$h=2^{-1}(m+1)^{-1}$, we conclude that
\[
  \|I_\rho-A_n\|_{\mathscr L(\sob,\R)}
  \ge\frac{I_\rho(f)}{\|f\|_{\sob}}
  \ge\frac{\rho_{\min}\Phi_0\,h^\alpha}{2(\rho_{\max}\Phi_\alpha)^{1/2}}
  =C_{\rho,\alpha}\,(m+1)^{-\alpha}.
\]
This concludes the proof.
\end{proof}

Since $m\le n$, Lemma~\ref{lem:fooling} implies in particular that any $n$-point
linear quadrature satisfies
\begin{equation}\label{eq:general-lower}
  \|I_\rho-A_n\|_{\mathscr L(\sob,\R)}\ge C_{\rho,\alpha}(n+1)^{-\alpha}.
\end{equation}
This recovers the general lower bound in \cite[Proposition~3.4]{SHK} for monotonic
Schwartz weights; see also \cite{DN24} for the Gaussian weight. Note that only the
positivity and continuity of $\rho$ on $[0,1]$ are used.

\subsection{Counting Gaussian nodes by Christoffel functions}\label{sec:count}

\begin{lemma}\label{lem:count}
Let $\rho$ be a weight and let $J\subset\R$ be a compact interval. Then the Gaussian nodes for $\rho$ satisfy
\[
  \#\{j\in\{1,\dots,n\}\;|\;x_{n,j}\in J\}\le\frac{\rho(\R)}{\min_{x\in J}\lambda_n(\rho;x)},
\]
where $\rho(\R):=\int_\R\rho(x)\dd x$.
\end{lemma}

\begin{proof}
By Lemma~\ref{lem:christoffel}(i), $\lambda_n(\rho;\cdot)=1/K_n(\cdot,\cdot)$ is continuous and positive, so that the minimum over $J$ exists and is positive. By Lemma~\ref{lem:gauss-weights},
\[
  \#\{j\;|\;x_{n,j}\in J\}\cdot\min_{x\in J}\lambda_n(\rho;x)\le\sum_{j:\,x_{n,j}\in J}\lambda_n(\rho;x_{n,j})=\sum_{j:\,x_{n,j}\in J}w_{n,j}\le\sum_{j=1}^nw_{n,j}=\rho(\R),
\]
which proves the claim.
\end{proof}

   We emphasize that the Gaussian quadrature weights enter only as a tool for counting the nodes. The lower bounds in Theorem~\ref{thm:main} concern the node set and hold for arbitrary quadrature weights $v_j$.

\subsection{Christoffel functions of the Lebesgue measure on an interval}\label{sec:interval}

In this and the following subsections, we bound $\lambda_n(\rho;x)$ from below for $|x|\le1$. We first need an elementary bound for the Lebesgue measure on an interval. For $L>0$, we write $\lambda_m([-L,L];x)$ for the Christoffel function of the measure $\mathbf 1_{[-L,L]}(t)\dd t$.

\begin{lemma}\label{lem:interval}
Let $m\in\N$ and $L>0$. Then, for $x\in(-L,L)$,
\[
  \lambda_m([-L,L];x)\ge\frac{\pi L}{2m}\Bigl(1-\frac{x^2}{L^2}\Bigr).
\]
In particular, $\lambda_m([-L,L];x)\ge3\pi L/(8m)$ for $|x|\le L/2$.
\end{lemma}

\begin{proof}
We first consider $L=1$. Since $\sqrt{1-t^2}\le1$ on $[-1,1]$, Lemma~\ref{lem:christoffel}(iii) shows that $\lambda_m([-1,1];y)$ is bounded below by the Christoffel function of the measure $\sqrt{1-t^2}\,\mathbf 1_{[-1,1]}(t)\dd t$.

The orthonormal polynomials with respect to the latter measure are $\sqrt{2/\pi}\,U_k$, where $U_k$ denotes the Chebyshev polynomial of the second kind, $U_k(\cos\theta)=\sin((k+1)\theta)/\sin\theta$ for $\theta\in(0,\pi)$. Indeed, the substitution $t=\cos\theta$ gives
\[
  \int_{-1}^1U_k(t)U_l(t)\sqrt{1-t^2}\dd t=\int_0^\pi\sin((k+1)\theta)\sin((l+1)\theta)\dd\theta=\frac\pi2\delta_{kl}.
\]
Since $|U_k(y)|\le(1-y^2)^{-1/2}$ for $y\in(-1,1)$, Lemma~\ref{lem:christoffel}(i) yields
\[
  \lambda_m([-1,1];y)\ge\Bigl(\frac2\pi\sum_{k=0}^{m-1}U_k(y)^2\Bigr)^{-1}\ge\frac{\pi(1-y^2)}{2m}.
\]

For general $L>0$, the map $q\mapsto\tilde q$ with $\tilde q(y):=q(Ly)$, $y\in\R$, is a bijection of $\Poly_{m-1}$ onto itself, and $q(x)=1$ if and only if $\tilde q(x/L)=1$. Moreover, $\int_{-L}^Lq(t)^2\dd t=L\int_{-1}^1\tilde q(y)^2\dd y$. Taking the infimum over all $q\in\Poly_{m-1}$ with $q(x)=1$, or equivalently over all $\tilde q\in\Poly_{m-1}$ with $\tilde q(x/L)=1$, we obtain $\lambda_m([-L,L];x)=L\,\lambda_m([-1,1];x/L)$. The claim now follows from the case $L=1$.
\end{proof}

\subsection{Freud weights with even exponents}\label{sec:even}

When $\beta=2r$ is an even integer, $|x|^\beta=x^{2r}$ is a polynomial. This allows us to construct an explicit polynomial that bounds $\rho_\beta$ from below on a large interval, and to apply Lemma~\ref{lem:christoffel}(iv).

\begin{proposition}\label{prop:even}
Let $r\in\N$, $\beta=2r$ and $\rho=\rho_\beta$. Then, for $n\ge n_r:=r\,4^r$ and $|x|\le1$,
\[
  \lambda_n(\rho;x)\ge\frac{\pi}{32}\,r^{-1/(2r)}\,n^{1/\beta-1}.
\]
Consequently, $\lambda_n(\rho;x)\ge c_r\,n^{1/\beta-1}$ for all $n\in\N$ and $|x|\le1$, where
\[
  c_r:=\frac{\pi}{32}\,r^{-1/(2r)}\,n_r^{1/(2r)-1}.
\]
\end{proposition}

\begin{proof}
Let $n\ge n_r$ and define
\[
  N:=2\Bigl\lceil\frac{n}{2r}\Bigr\rceil,\qquad L:=N^{1/(2r)},\qquad G(t):=\Bigl(1-\frac{t^{2r}}{N}\Bigr)^{N/2}.
\]
Then $N$ is an even integer with $n/r\le N<n/r+2$, and $G\in\Poly_{rN}$. Since $n\ge r4^r$, we have
\[
  L\ge(n/r)^{1/(2r)}\ge2,\qquad rN<n+2r\le2n .
\]

For $|t|\le L$ we have $s:=t^{2r}/N\in[0,1]$, and since $0\le1-s\le\mathrm{e}^{-s}$,
\[
  G(t)^2=(1-s)^N\le\mathrm{e}^{-sN}=\rho(t).
\]
Hence $G^2\mathbf 1_{[-L,L]}\le\rho$ pointwise, and the assumption of Lemma~\ref{lem:christoffel}(iv) holds with $\nu=\rho(t)\dd t$, $\nu'=\mathbf 1_{[-L,L]}(t)\dd t$ and $k=rN$.

Moreover, for $|x|\le1$ we have $s\le1/N$, so that $G(x)^2\ge(1-1/N)^N\ge1/4$, since $N\ge2$ and $N\mapsto(1-1/N)^N$ is increasing.

Since $|x|\le1\le L/2$, Lemma~\ref{lem:christoffel}(iv) and Lemma~\ref{lem:interval} yield
\[
  \lambda_n(\rho;x)\ge G(x)^2\,\lambda_{n+rN}([-L,L];x)\ge\frac14\cdot\frac{3\pi L}{8(n+rN)}\ge\frac{\pi L}{32\,n}\ge\frac{\pi}{32}\Bigl(\frac nr\Bigr)^{1/(2r)}\frac1n ,
\]
where we used $n+rN<3n$. This proves the first claim.

For $1\le n<n_r$, Lemma~\ref{lem:christoffel}(ii) gives $\lambda_n(\rho;x)\ge\lambda_{n_r}(\rho;x)\ge\frac{\pi}{32}r^{-1/(2r)}n_r^{1/(2r)-1}$, and the second claim follows since $n^{1/\beta-1}\le1$.
\end{proof}

\subsection{General Freud weights}\label{sec:freud}

For non-even $\beta$, the above construction is not available. We instead rely on two standard results on weighted polynomials. Following the convention in the literature, we write
\[
  Q(x):=\tfrac12|x|^\beta,\qquad W(x):=\mathrm{e}^{-Q(x)},\qquad\text{so that}\qquad\rho_\beta=W^2.
\]
Since $\beta>1$, $Q$ is continuously differentiable with $Q'(x)=\frac\beta2\,\sign(x)|x|^{\beta-1}$.

The Mhaskar--Rakhmanov--Saff number $a_n$ is defined as the positive solution of
\[
  n=\frac2\pi\int_0^1\frac{a_nu\,Q'(a_nu)}{\sqrt{1-u^2}}\dd u .
\]
For our $Q$, we have $a_nuQ'(a_nu)=\frac\beta2a_n^\beta u^\beta$, and the substitution $s=u^2$ gives
\[
  \int_0^1\frac{u^\beta}{\sqrt{1-u^2}}\dd u=\frac{\sqrt\pi\,\Gamma(\frac{\beta+1}2)}{2\,\Gamma(\frac\beta2+1)}.
\]
Hence
\begin{equation}\label{eq:mrs}
  a_n=\Bigl(\frac{n}{\theta_\beta}\Bigr)^{1/\beta},\qquad\theta_\beta:=\frac{\beta\,\Gamma(\frac{\beta+1}2)}{2\sqrt\pi\,\Gamma(\frac\beta2+1)},
\end{equation}
so that $a_n^\beta=n/\theta_\beta$ and $a_n/n=\theta_\beta^{-1/\beta}n^{1/\beta-1}$. For $\beta=2$, \eqref{eq:mrs} gives $a_n=\sqrt{2n}$.

We recall the following two results.

\begin{lemma}[Restricted range inequality \cite{MS85}; see also \cite{LL01,ST97}]\label{lem:range}
For any $n\in\N$ and $P\in\Poly_n$,
\[
  \|PW\|_{L^\infty(\R)}=\|PW\|_{L^\infty([-a_n,a_n])}.
\]
\end{lemma}

\begin{lemma}[Markov--Bernstein inequality \cite{LL90}; see also \cite{LL01,Lub07}]\label{lem:bernstein}
There exists a constant $C_\beta\ge1$, depending only on $\beta$, such that, for any $n\in\N$ and $P\in\Poly_n$,
\[
  \|P'W\|_{L^\infty(\R)}\le C_\beta\,\frac{n}{a_n}\,\|PW\|_{L^\infty(\R)}.
\]
\end{lemma}

\begin{proposition}\label{prop:freud}
Let $\beta>1$ and $\rho=\rho_\beta$. Let $C_\beta$ be the constant in Lemma~\ref{lem:bernstein} and set $D_\beta:=C_\beta+\beta\,2^{\beta-2}/\theta_\beta$. Then, for all $n\in\N$ and $x\in\R$,
\[
  \lambda_n(\rho;x)\ge\frac{1}{4D_\beta}\,\frac{a_n}{n}\,\rho(x).
\]
In particular, $\lambda_n(\rho;x)\ge(4\mathrm e D_\beta)^{-1}\theta_\beta^{-1/\beta}\,n^{1/\beta-1}$ for $|x|\le1$.
\end{proposition}

\begin{proof}
Let $q\in\Poly_{n-1}$ with $q(x)=1$. Define $F:=qW$ and $M:=\|F\|_{L^\infty(\R)}$. Then $F$ is continuously differentiable and $M\ge|F(x)|=W(x)>0$. By Lemma~\ref{lem:range} with $P=q\in\Poly_n$, $M=\|F\|_{L^\infty([-a_n,a_n])}$. Hence, by continuity, there exists $x^*\in[-a_n,a_n]$ with $|F(x^*)|=M$.

We now bound the derivative $F'=q'W-Q'qW$. By Lemma~\ref{lem:bernstein}, $|q'(t)W(t)|\le C_\beta(n/a_n)M$ for all $t\in\R$. Moreover, for $|t|\le2a_n$,
\[
  |Q'(t)|\le\frac\beta2(2a_n)^{\beta-1}=\beta\,2^{\beta-2}\,\frac{a_n^\beta}{a_n}=\frac{\beta\,2^{\beta-2}}{\theta_\beta}\,\frac{n}{a_n}.
\]
Therefore, noting $D_\beta\ge1$,
\[
  |F'(t)|\le D_\beta\,\frac{n}{a_n}\,M\qquad\text{for all }|t|\le2a_n .
\]

Set $\delta:=a_n/(2D_\beta n)\le a_n$, so that $[x^*-\delta,x^*+\delta]\subset[-2a_n,2a_n]$. By the mean value theorem, for $|t-x^*|\le\delta$,
\[
  |F(t)|\ge|F(x^*)|-D_\beta\frac{n}{a_n}M|t-x^*|\ge\frac M2 .
\]
Consequently,
\[
  \int_\R q(t)^2\rho(t)\dd t=\int_\R F(t)^2\dd t\ge\int_{x^*-\delta}^{x^*+\delta}\frac{M^2}{4}\dd t=\frac{M^2a_n}{4D_\beta n}\ge\frac{a_n}{4D_\beta n}\,W(x)^2=\frac{a_n}{4D_\beta n}\,\rho(x).
\]
Taking the infimum over all $q\in\Poly_{n-1}$ with $q(x)=1$ proves the first claim. The second claim follows from $\rho(x)\ge\mathrm{e}^{-1}$ for $|x|\le1$ and \eqref{eq:mrs}.
\end{proof}

\begin{remark}\label{rem:two-sided}
Two-sided estimates of the form $c_\beta(a_n/n)\rho_\beta(x)\le\lambda_n(\rho_\beta;x)\le C'_\beta(a_n/n)\rho_\beta(x)$, valid uniformly for $|x|\le a_n/2$ with constants depending only on $\beta$, are known \cite{LL92,LL01}, so that Proposition~\ref{prop:freud} may also be replaced by a direct citation. We have included the above derivation since it shows transparently where the scale $a_n/n$ originates.
\end{remark}

\subsection{The logistic weight}\label{sec:logistic}

For the logistic weight, the three-term recurrence coefficients of the orthogonal polynomials are explicitly available. This allows us to compute the Christoffel--Darboux kernel at the origin in closed form.

Let $(P_n)_{n\ge0}$ be the monic orthogonal polynomials with respect to $\rho_{\mathrm L}$. Since $\rho_{\mathrm L}$ is even, $P_n(-x)=(-1)^nP_n(x)$, and $(P_n)$ satisfies the three-term recurrence
\[
  xP_n(x)=P_{n+1}(x)+b_nP_{n-1}(x),\qquad P_{-1}:=0,\quad P_0:=1,
\]
with $b_n>0$. Moreover, $\|P_n\|_{L^2_\rho(\R)}^2=b_1\cdots b_n$, and since $\rho_{\mathrm L}(\R)=1$, the orthonormal polynomials are $p_n=P_n/\|P_n\|_{L^2_\rho(\R)}$ with $p_0=1$.

\begin{lemma}\label{lem:recurrence}
For the logistic weight $\rho_{\mathrm L}$,
\[
  b_n=\frac{\pi^2n^4}{4n^2-1}\qquad(n\in\N).
\]
\end{lemma}

\begin{proof}
By the reflection formula, $|\Gamma(\frac12+\mathrm it)|^2=\pi/\cosh(\pi t)$ for $t\in\R$. Hence, under the change of variables $x=2\pi t$,
\[
  2\pi\,\rho_{\mathrm L}(2\pi t)=\frac{\pi}{2\cosh^2(\pi t)}=\frac{1}{2\pi}\bigl|\Gamma(\tfrac12+\mathrm it)\bigr|^4 .
\]
The right-hand side is the weight of the continuous Hahn polynomials with parameters $a=b=c=d=\frac12$; see \cite[Section~9.4]{KLS10} or \cite[(1.4.2)]{KS98}.

The monic continuous Hahn polynomials $\tilde p_n$ (denoted by $p_n$ in \cite{KS98}) satisfy \cite[Section~9.4]{KLS10}, \cite[(1.4.3)--(1.4.4)]{KS98}
\[
  t\,\tilde p_n(t)=\tilde p_{n+1}(t)+\mathrm i(A_n+C_n+a)\tilde p_n(t)-A_{n-1}C_n\tilde p_{n-1}(t),
\]
where, with $s:=a+b+c+d$,
\[
  A_n=-\frac{(n+s-1)(n+a+c)(n+a+d)}{(2n+s-1)(2n+s)},\qquad C_n=\frac{n(n+b+c-1)(n+b+d-1)}{(2n+s-2)(2n+s-1)}.
\]
For $a=b=c=d=\frac12$ we have $s=2$, and
\[
  A_n=-\frac{(n+1)^2}{2(2n+1)},\qquad C_n=\frac{n^2}{2(2n+1)}.
\]
Thus $A_n+C_n+\frac12=0$, and $-A_{n-1}C_n=n^4/(4(4n^2-1))$.

Since $P_n(x)=(2\pi)^n\tilde p_n(x/(2\pi))$, the recurrence for $\tilde p_n$ transforms into that for $P_n$ with $b_n=(2\pi)^2\cdot n^4/(4(4n^2-1))$, which proves the claim.
\end{proof}

As a consistency check, we compare $b_1$ and $b_2$ with the moments
$m_k:=\int_\R x^k\rho_{\mathrm L}(x)\dd x$, $k\in\N\cup\{0\}$, where $m_0=1$. The
recurrence gives $P_1(x)=x$ and $P_2(x)=x^2-b_1$. Together with
$\|P_n\|_{L^2_\rho(\R)}^2=b_1\cdots b_n$, this yields
\[
  b_1=\int_\R x^2\rho_{\mathrm L}(x)\dd x=m_2,
  \qquad
  b_1b_2=\int_\R\bigl(x^2-b_1\bigr)^2\rho_{\mathrm L}(x)\dd x=m_4-2b_1m_2+b_1^2=m_4-m_2^2,
\]
that is, $b_1=m_2$ and $b_2=(m_4-m_2^2)/m_2$. Since
\[
  \int_\R\mathrm e^{sx}\rho_{\mathrm L}(x)\dd x=\frac{\pi s}{\sin(\pi s)}
  =1+\frac{\pi^2s^2}{6}+\frac{7\pi^4s^4}{360}+O(s^6),\qquad |s|<1,
\]
we have $m_2=2!\cdot\pi^2/6=\pi^2/3$ and $m_4=4!\cdot7\pi^4/360=7\pi^4/15$. Hence
$b_1=\pi^2/3$ and $b_2=16\pi^2/15$, in agreement with Lemma~\ref{lem:recurrence}.

In what follows, we write $\gamma_m:=4^{-m}\binom{2m}{m}=\prod_{k=1}^m\frac{2k-1}{2k}$ for $m\in\N$, and $\gamma_0:=1$.

\begin{lemma}\label{lem:origin}
For the logistic weight and $m\in\N\cup\{0\}$, the orthonormal polynomials satisfy
\[
  p_{2m+1}(0)=0,\qquad p_{2m}(0)^2=(4m+1)\,\gamma_m^4 .
\]
Consequently, for $n\in\N$,
\[
  1+\frac14\log\frac n2\le K_n(0,0)=\sum_{m=0}^{\lfloor(n-1)/2\rfloor}(4m+1)\gamma_m^4\le1+\frac5{\pi^2}(1+\log n),
\]
and $\lim_{n\to\infty}K_n(0,0)/\log n=4/\pi^2$.
\end{lemma}

\begin{proof}
By symmetry, $P_{2m+1}(0)=0$. Evaluating the recurrence at $x=0$ gives $P_{2m}(0)=-b_{2m-1}P_{2m-2}(0)$. Hence
\[
  p_{2m}(0)^2=\frac{b_{2m-1}^2P_{2m-2}(0)^2}{b_1\cdots b_{2m}}=p_{2m-2}(0)^2\,\frac{b_{2m-1}}{b_{2m}} .
\]
Using $4n^2-1=(2n-1)(2n+1)$ and Lemma~\ref{lem:recurrence}, we obtain
\[
  \frac{b_{2m-1}}{b_{2m}}=\Bigl(\frac{2m-1}{2m}\Bigr)^4\frac{4m+1}{4m-3}.
\]
Starting from $p_0(0)^2=1$ and multiplying these ratios, the second factors telescope to $4m+1$, which proves the first claim.

For the bounds on $K_n(0,0)$, note that $\gamma_m\sqrt m$ is increasing in $m\in\N$, since
\[
  \frac{\gamma_{m+1}^2(m+1)}{\gamma_m^2m}=\frac{(2m+1)^2}{4m(m+1)}>1 .
\]
Moreover, $\gamma_1=1/2$, and Wallis' formula
$\prod_{k=1}^\infty\frac{(2k)^2}{(2k-1)(2k+1)}=\frac\pi2$ gives
\[
  (2m+1)\gamma_m^2=\prod_{k=1}^m\frac{(2k-1)(2k+1)}{(2k)^2}\to\frac2\pi ,
\]
so that $\gamma_m\sqrt m\to\pi^{-1/2}$ as $m\to\infty$. Hence $(2\sqrt m)^{-1}\le\gamma_m\le(\pi m)^{-1/2}$ for $m\in\N$, and therefore
\[
  \frac{1}{4m}=\frac{4m}{16m^2}\le(4m+1)\gamma_m^4\le\frac{4m+1}{\pi^2m^2}\le\frac5{\pi^2m}.
\]
With $M:=\lfloor(n-1)/2\rfloor$, we obtain
\[
  1+\frac14\sum_{m=1}^{M}\frac1m\le K_n(0,0)\le1+\frac5{\pi^2}\sum_{m=1}^{M}\frac1m .
\]
The bounds follow from $\log(M+1)\le\sum_{m=1}^M1/m\le1+\log n$ and $M+1\ge n/2$. Finally, since $m\gamma_m^2\to1/\pi$, we have $M(4M+1)\gamma_M^4\to4/\pi^2$, which is
the ratio of the $M$th terms of $\sum_{m=1}^M(4m+1)\gamma_m^4$ and
$\sum_{m=1}^M1/m$. Since the latter sum diverges, the Stolz--Ces\`aro theorem shows
that the ratio of the two sums also tends to $4/\pi^2$. The limit follows from
$\sum_{m=1}^M1/m\sim\log M\sim\log n$. 
\end{proof}

   By Lemma~\ref{lem:origin}, $\lambda_n(\rho_{\mathrm L};0)\log n\to\pi^2/4$ as $n\to\infty$, so that the Gaussian quadrature weights near the origin decay only logarithmically. Since $w_{n,j}$ is approximately $\rho_{\mathrm L}(x_{n,j})$ times the local node spacing, the spacing of the Gaussian nodes near the origin is heuristically close to $\pi^2/\log n$. To extend the bound from the origin to $[0,1]$, we use a translation argument.

\begin{lemma}\label{lem:shift}
For all $n\in\N$ and $x,t\in\R$, it holds that $\lambda_n(\rho_{\mathrm L};x+t)\ge\mathrm{e}^{-|t|}\lambda_n(\rho_{\mathrm L};x)$.
\end{lemma}

\begin{proof}
By the addition formula, $\cosh\frac{y+t}2\le\cosh\frac y2\,(\cosh\frac t2+|\sinh\frac t2|)=\mathrm{e}^{|t|/2}\cosh\frac y2$. Hence $\rho_{\mathrm L}(y+t)\ge\mathrm{e}^{-|t|}\rho_{\mathrm L}(y)$ for all $y\in\R$.

Let $q\in\Poly_{n-1}$ with $q(x+t)=1$, and set $\tilde q(u):=q(u+t)$, $u\in\R$. Then $\tilde q\in\Poly_{n-1}$ and $\tilde q(x)=1$, so that
\[
  \int_\R q(y)^2\rho_{\mathrm L}(y)\dd y
  =\int_\R\tilde q(u)^2\rho_{\mathrm L}(u+t)\dd u
  \ge\mathrm{e}^{-|t|}\int_\R\tilde q(u)^2\rho_{\mathrm L}(u)\dd u
  \ge\mathrm{e}^{-|t|}\lambda_n(\rho_{\mathrm L};x).
\]
Taking the infimum over all $q\in\Poly_{n-1}$ with $q(x+t)=1$ completes the proof.
\end{proof}

\begin{proposition}\label{prop:logistic}
For all $n\in\N$ and $|x|\le1$,
\[
  \frac{1}{\mathrm e\,\bigl(1+\frac5{\pi^2}(1+\log n)\bigr)}\le\lambda_n(\rho_{\mathrm L};x)\le\frac{\mathrm e}{1+\frac14\log\frac n2} .
\]
\end{proposition}

\begin{proof}
Lemma~\ref{lem:shift} with base point $0$ and shift $x$ gives $\lambda_n(\rho_{\mathrm L};x)\ge\mathrm{e}^{-|x|}\lambda_n(\rho_{\mathrm L};0)$. Lemma~\ref{lem:shift} with base point $x$ and shift $-x$ gives $\lambda_n(\rho_{\mathrm L};0)\ge\mathrm{e}^{-|x|}\lambda_n(\rho_{\mathrm L};x)$. Since $\lambda_n(\rho_{\mathrm L};0)=1/K_n(0,0)$ by Lemma~\ref{lem:christoffel}(i), both claims follow from Lemma~\ref{lem:origin}.
\end{proof}

\subsection{Proof of Theorem~\ref{thm:main}}\label{sec:proof-main}

\begin{proof}[Proof of Theorem~\ref{thm:main}]
Let $m_n$ denote the number of Gaussian nodes in $(0,1)$.

(i) Let $\rho=\rho_\beta$. By Proposition~\ref{prop:even} (when $\beta$ is an even integer) or Proposition~\ref{prop:freud} (for general $\beta>1$), there exists $c_\beta>0$, depending only on $\beta$, such that $\lambda_n(\rho;x)\ge c_\beta\,n^{1/\beta-1}$ for all $n\in\N$ and $|x|\le1$. Lemma~\ref{lem:count} with $J=[0,1]$ then gives
\[
  m_n\le\frac{2\Gamma(1+1/\beta)}{c_\beta}\,n^{1-1/\beta}=:C'_\beta\,n^{1-1/\beta}.
\]
Applying Lemma~\ref{lem:fooling} with $\rho_{\max}=1$ and $\rho_{\min}=\mathrm{e}^{-1}$, and using $n^{1-1/\beta}\ge1$, we obtain
\[
  \|I_\rho-A_n\|_{\mathscr L(\sob,\R)}\ge C_{\rho,\alpha}(C'_\beta n^{1-1/\beta}+1)^{-\alpha}\ge C_{\rho,\alpha}(C'_\beta+1)^{-\alpha}\,n^{-\alpha(1-1/\beta)} .
\]

(ii) Let $\rho=\rho_{\mathrm L}$ and $n\ge2$. By Proposition~\ref{prop:logistic} and Lemma~\ref{lem:count} with $J=[0,1]$ and $\rho_{\mathrm L}(\R)=1$,
\[
  m_n+1\le1+\mathrm e\Bigl(1+\frac5{\pi^2}(1+\log n)\Bigr)\le C_{\mathrm L}\log n ,\qquad C_{\mathrm L}:=\frac{1+\mathrm e(1+5/\pi^2)}{\log2}+\frac{5\mathrm e}{\pi^2},
\]
where we used $1\le\log n/\log2$ for $n\ge2$. Lemma~\ref{lem:fooling} with $\rho_{\max}=\rho_{\mathrm L}(0)=1/4$ and $\rho_{\min}=\rho_{\mathrm L}(1)=1/(4\cosh^2(1/2))$ then gives
\[
  \|I_\rho-A_n\|_{\mathscr L(\sob,\R)}\ge C_{\rho,\alpha}\,C_{\mathrm L}^{-\alpha}(\log n)^{-\alpha}.
\]
This concludes the proof.
\end{proof}

\section{Optimal quadrature by the M\"obius-transformed trapezoidal rule}\label{sec:mtr}

In this section, we show that, for the logistic weight and for the Freud weights with
even integer exponents, the rates in Theorem~\ref{thm:main} are far from optimal, i.e., the optimal rate
follows directly from \cite{SHK}. For the Freud weights with non-even exponents,
$\rho_\beta$ is not smooth at the origin, and Theorem~\ref{thm:shk} does not apply. However, 
the optimal rate $n^{-\alpha}$ is known from \cite[Corollary~3.3]{Dun25}.

For $c>0$, let
\[
  \phi_c(\theta):=-c\cot(\theta/2),\qquad\phi_c'(\theta)=\frac{c}{2\sin^2(\theta/2)},\qquad\theta\in(0,2\pi),
\]
which maps $(0,2\pi)$ increasingly onto $\R$. For a weight $\rho$ and $n\in\N$, the M\"obius-transformed trapezoidal rule \cite{SHK} is defined by
\[
  Q^{\mathrm{MT}}_{\rho,n}(f):=\frac{2\pi}{n}\sum_{j=1}^{n-1}f(\phi_c(\theta_j))\,\rho(\phi_c(\theta_j))\,\phi_c'(\theta_j),\qquad\theta_j:=\frac{2\pi j}{n}.
\]
This corresponds to the trapezoidal rule for the $2\pi$-periodic function $((f\rho)\circ\phi_c)\,\phi_c'$, where the term at $\theta=0\equiv2\pi$ is set to zero. 
Note that $Q^{\mathrm{MT}}_{\rho,n}$ uses $n-1$ function values and only requires point evaluations of $\rho$. Throughout this section, $c>0$ is fixed, and all constants may depend on $c$ without further mention.

\begin{theorem}[{\cite[Theorem~3.1]{SHK}}]\label{thm:shk}
Let $\alpha\in\N$ and $\rho\in\mathcal S^{\mathrm{mon}}_+$. Then there exists a constant
$C^{\mathrm{MT}}_{\rho,\alpha}>0$, depending only on $\rho$, $\alpha$ and $c$, such that for
all $n\in\N$ and $f\in W^{\alpha,2}_\rho(\R)$,
\[
  |I_\rho(f)-Q^{\mathrm{MT}}_{\rho,n}(f)|\le C^{\mathrm{MT}}_{\rho,\alpha}\,n^{-\alpha}\,
  \|f\|_{W^{\alpha,2}_\rho(\R)}.
\]
\end{theorem}

\begin{lemma}\label{lem:schwartz}
The logistic weight $\rho_{\mathrm L}$ and the Freud weights $\rho_{2r}$, $r\in\N$, belong to $\mathcal S^{\mathrm{mon}}_+$.
\end{lemma}

\begin{proof}
For $\rho_{2r}$, every derivative is a polynomial multiple of $\mathrm{e}^{-x^{2r}}$, and $x\rho_{2r}'(x)=-2r\,x^{2r}\rho_{2r}(x)\le0$.

For $\rho_{\mathrm L}=\frac14\sech^2(x/2)$, set $s=\sech(x/2)$ and $\sigma=\tanh(x/2)$, so that $s'=-\frac12s\sigma$ and $\sigma'=\frac12s^2$. By induction, every derivative of $\rho_{\mathrm L}$ is a linear combination of monomials $s^k\sigma^l$ with $k\ge2$. Since $|\sigma|\le1$ and $s\le2\mathrm{e}^{-|x|/2}$, every derivative of $\rho_{\mathrm L}$ is bounded by a constant multiple of $\mathrm{e}^{-|x|}$, so that $\rho_{\mathrm L}\in\mathcal S$. Moreover, $x\rho_{\mathrm L}'(x)=-\frac14xs^2\sigma\le0$.
\end{proof}

Combining Theorem~\ref{thm:shk}, Lemma~\ref{lem:schwartz}, the general lower bound \eqref{eq:general-lower} and Theorem~\ref{thm:main}, we obtain the following.

\begin{corollary}\label{cor:mtr-smooth}
Let $\alpha\in\N$, and let $\rho=\rho_{\mathrm L}$ or $\rho=\rho_{2r}$ with $r\in\N$. Then, for all $n\ge2$,
\[
  C_{\rho,\alpha}\,n^{-\alpha}\le\|I_\rho-Q^{\mathrm{MT}}_{\rho,n}\|_{\mathscr L(\sob,\R)}\le C^{\mathrm{MT}}_{\rho,\alpha}\,n^{-\alpha},
\]
where $C_{\rho,\alpha}$ is the constant in Lemma~\ref{lem:fooling}. Moreover,
\[
  \frac{\|I_\rho-Q^{\mathrm G}_{\rho,n}\|_{\mathscr L(\sob,\R)}}{\|I_\rho-Q^{\mathrm{MT}}_{\rho,n}\|_{\mathscr L(\sob,\R)}}\ge
  \begin{cases}
    \dfrac{c_{\alpha,2r}}{C^{\mathrm{MT}}_{\rho,\alpha}}\,n^{\alpha/(2r)}, & \rho=\rho_{2r},\\[10pt]
    \dfrac{c_{\alpha}}{C^{\mathrm{MT}}_{\rho,\alpha}}\,n^{\alpha}(\log n)^{-\alpha}, & \rho=\rho_{\mathrm L},
  \end{cases}
\]
with the constants $c_{\alpha,2r}$ and $c_\alpha$ from Theorem~\ref{thm:main}.
\end{corollary}

\section{Concluding remarks}\label{sec:concluding}

In this paper, we proved lower bounds on the worst-case error of Gaussian quadrature for Freud weights and for the logistic weight in weighted Sobolev spaces, and showed that the M\"obius-transformed trapezoidal rule, with a simple modification if necessary, attains the optimal rate $n^{-\alpha}$ for all these weights.

The rates in Theorem~\ref{thm:main} have a simple interpretation. By \eqref{eq:mrs}, we have $n^{1-1/\beta}=\theta_\beta^{1/\beta}\,n/a_n$. Since the $n$ Gaussian nodes for $\rho_\beta$ are spread over $[-a_n,a_n]$ with $a_n=(n/\theta_\beta)^{1/\beta}$, a fixed interval contains at most a constant multiple of $n/a_n$ nodes, and Lemma~\ref{lem:fooling} converts this node count into the rate $n^{-\alpha(1-1/\beta)}$. For $\beta=2$ this is the rate $n^{-\alpha/2}$ of \cite{KSG23}; as $\beta\to\infty$ the exponent approaches the optimal exponent $\alpha$, while it tends to zero as $\beta\to1$. The logistic weight decays exactly as the borderline Freud weight $\mathrm{e}^{-|x|}$, in which case the number of nodes in a fixed interval is at most a constant multiple of $\log n$. This is the origin of the slow convergence of Gauss--Logistic quadrature observed numerically in \cite{SHK}.

We also expect our results to extend to multivariate integration. For the Gaussian weight, it was shown in \cite{KSG_sparse} that, in Sobolev spaces of dominating mixed smoothness, several quasi-Monte Carlo methods with a change of variables attain the optimal rate, whereas the dominant rate of sparse-grid quadrature based on the Gauss--Hermite rule is only half of the optimal one. We expect the same picture for the weights considered here: quasi-Monte Carlo methods combined with the M\"obius transformation should attain the optimal rate (see also \cite[Section~6]{SHK}), whereas sparse-grid quadrature based on Gaussian rules should not. Since the lower bound in \cite{KSG_sparse} is based on a univariate fooling function, we expect that Lemma~\ref{lem:fooling} and the node counts in Section~\ref{sec:main} can be used in a similar way.

We mention some directions for future work. First, it remains open whether the
M\"obius-transformed trapezoidal rule, possibly combined with a change of variables
that is flat at the origin, attains the optimal rate $n^{-\alpha}$ for Freud weights
with non-even exponents. Second, it would be natural to consider Sobolev spaces based
on $L^p$ with $p\ne2$. For Freud weights, the optimal rate in such spaces is again
$n^{-\alpha}$ for $1<p\le\infty$ \cite[Corollary~3.3]{Dun25}, whereas for $p=1$ it is
$n^{-\alpha(1-1/\beta)}$ and is attained by truncated Gaussian rules \cite{Dun23}.
It would be interesting to determine the convergence rate of Gaussian quadrature in
these spaces. Third, apart from the Gaussian weight, for which matching upper bounds
are known \cite{KSG23}, our lower bounds for Gaussian quadrature are not accompanied
by upper bounds. Establishing such upper bounds for general Freud weights and for the
logistic weight would show that the lower bounds in Theorem~\ref{thm:main} are sharp.

\section*{Declaration of use of AI tools}
The author used Claude Opus 5.5 (Anthropic) as an interactive assistant during the development and preparation of this work. Its use includes developing and revising proof arguments (Sections~\ref{sec:interval}--\ref{sec:logistic}), checking calculations and consistency, and editing the exposition and
LaTeX source. 
The content was independently reviewed and verified by the author. The author assumes full responsibility for all content.

\bibliographystyle{plain}
\bibliography{ref_subopt}

\begin{thebibliography}{10}

\bibitem{Bak67}
Nikolai~S. Bakhvalov.
\newblock On the optimal speed of integrating analytic functions.
\newblock {\em USSR Comput. Math. Math. Phys.}, 7(5):63--75, 1967.

\bibitem{Cur98}
Francisco Curbera.
\newblock Optimal integration of {L}ipschitz functions with a {G}aussian
  weight.
\newblock {\em J. Complexity}, 14(1):122--149, 1998.

\bibitem{DHP15}
Josef Dick, Aicke Hinrichs, and Friedrich Pillichshammer.
\newblock Proof techniques in quasi-{M}onte {C}arlo theory.
\newblock {\em J. Complexity}, 31(3):327--371, 2015.

\bibitem{DIL18}
Josef Dick, Christian Irrgeher, Gunther Leobacher, and Friedrich
  Pillichshammer.
\newblock On the optimal order of integration in {H}ermite spaces with finite
  smoothness.
\newblock {\em SIAM J. Numer. Anal.}, 56(2):684--707, 2018.

\bibitem{Dun23}
Dinh D{\~u}ng.
\newblock Numerical weighted integration of functions having mixed smoothness.
\newblock {\em J. Complexity}, 78:101757, 2023.

\bibitem{Dun25}
Dinh D{\~u}ng.
\newblock Sampling reconstruction and integration of functions on
  {$\mathbb{R}^d$} endowed with a measure.
\newblock VIASM Preprint ViAsM25.15, 2025.

\bibitem{DN24}
Dinh D{\~u}ng and Van~Kien Nguyen.
\newblock Optimal numerical integration and approximation of functions on
  {$\mathbb{R}^d$} equipped with {G}aussian measure.
\newblock {\em IMA J. Numer. Anal.}, 44(2):1242--1267, 2024.

\bibitem{Fre71}
G{\'e}za Freud.
\newblock {\em Orthogonal Polynomials}.
\newblock Pergamon Press, Oxford, 1971.

\bibitem{Gau04}
Walter Gautschi.
\newblock {\em Orthogonal Polynomials: Computation and Approximation}.
\newblock Numerical Mathematics and Scientific Computation. Oxford University
  Press, New York, 2004.

\bibitem{GKT24}
Takashi Goda, Yoshihito Kazashi, and Ken'ichiro Tanaka.
\newblock How sharp are error bounds? --{L}ower bounds on quadrature worst-case
  errors for analytic functions--.
\newblock {\em SIAM J. Numer. Anal.}, 62(5):2370--2392, 2024.

\bibitem{HS2026}
Nuutti Hyv{\"o}nen and Yuya Suzuki.
\newblock {M}{\"o}bius-transformed trapezoidal rule for polynomial weights.
\newblock arXiv:2604.27758 [math.NA], 2026.

\bibitem{KSG23}
Yoshihito Kazashi, Yuya Suzuki, and Takashi Goda.
\newblock Suboptimality of {G}auss--{H}ermite quadrature and optimality of the
  trapezoidal rule for functions with finite smoothness.
\newblock {\em SIAM J. Numer. Anal.}, 61(3):1426--1448, 2023.

\bibitem{KSG_sparse}
Yoshihito Kazashi, Yuya Suzuki, and Takashi Goda.
\newblock Optimality of quasi-{M}onte {C}arlo methods and suboptimality of the
  sparse-grid {G}auss--{H}ermite rule in {G}aussian {S}obolev spaces.
\newblock {\em SIAM J. Numer. Anal.}, 64(5):2123--2150, 2026.

\bibitem{KLS10}
Roelof Koekoek, Peter~A. Lesky, and Ren{\'e}~F. Swarttouw.
\newblock {\em Hypergeometric Orthogonal Polynomials and Their
  {$q$}-Analogues}.
\newblock Springer Monographs in Mathematics. Springer-Verlag, Berlin, 2010.

\bibitem{KS98}
Roelof Koekoek and Ren{\'e}~F. Swarttouw.
\newblock The {A}skey-scheme of hypergeometric orthogonal polynomials and its
  $q$-analogue.
\newblock Report 98-17, Delft University of Technology, Faculty of Information
  Technology and Systems, 1998.
\newblock Available at
  \url{https://tudelft.roelofkoekoek.nl/documents/as98.pdf}.

\bibitem{LL90}
Eli Levin and Doron~S. Lubinsky.
\newblock {$L_\infty$} {M}arkov and {B}ernstein inequalities for {F}reud
  weights.
\newblock {\em SIAM J. Math. Anal.}, 21:1065--1082, 1990.

\bibitem{LL92}
Eli Levin and Doron~S. Lubinsky.
\newblock Christoffel functions, orthogonal polynomials, and {N}evai's
  conjecture for {F}reud weights.
\newblock {\em Constr. Approx.}, 8:463--535, 1992.

\bibitem{LL01}
Eli Levin and Doron~S. Lubinsky.
\newblock {\em Orthogonal Polynomials for Exponential Weights}.
\newblock CMS Books in Mathematics. Springer-Verlag, New York, 2001.

\bibitem{Lub07}
Doron~S. Lubinsky.
\newblock A survey of weighted polynomial approximation with exponential
  weights.
\newblock {\em Surv. Approx. Theory}, 3:1--105, 2007.

\bibitem{MM08}
Giuseppe Mastroianni and Gradimir~V. Milovanovi{\'c}.
\newblock {\em Interpolation Processes: Basic Theory and Applications}.
\newblock Springer Monographs in Mathematics. Springer, Berlin, Heidelberg,
  2008.

\bibitem{MS85}
Hrushikesh~N. Mhaskar and Edward~B. Saff.
\newblock Where does the sup norm of a weighted polynomial live? ({A}
  generalization of incomplete polynomials).
\newblock {\em Constr. Approx.}, 1:71--91, 1985.

\bibitem{Nev86}
Paul Nevai.
\newblock G{\'e}za {F}reud, orthogonal polynomials and {C}hristoffel functions.
  {A} case study.
\newblock {\em J. Approx. Theory}, 48:3--167, 1986.

\bibitem{Pet98}
Knut Petras.
\newblock Gaussian versus optimal integration of analytic functions.
\newblock {\em Constr. Approx.}, 14:231--245, 1998.

\bibitem{ST97}
Edward~B. Saff and Vilmos Totik.
\newblock {\em Logarithmic Potentials with External Fields}.
\newblock Grundlehren der mathematischen Wissenschaften. Springer-Verlag,
  Berlin, 1997.

\bibitem{SHK}
Yuya Suzuki, Nuutti Hyv{\"o}nen, and Toni Karvonen.
\newblock {M}{\"o}bius-transformed trapezoidal rule.
\newblock {\em Math. Comp.}, 95:1491--1515, 2026.

\bibitem{Sze}
G{\'a}bor Szeg{\H{o}}.
\newblock {\em Orthogonal Polynomials}, volume~23 of {\em American Mathematical
  Society Colloquium Publications}.
\newblock American Mathematical Society, Providence, RI, fourth edition, 1975.

\bibitem{Tem92}
Vladimir~N. Temlyakov.
\newblock On a way of obtaining lower estimates for the errors of quadrature
  formulas.
\newblock {\em Math. USSR-Sb.}, 71(1):247--257, 1992.

\bibitem{Tre22}
Lloyd~N. Trefethen.
\newblock Exactness of quadrature formulas.
\newblock {\em SIAM Rev.}, 64(1):132--150, 2022.

\end{thebibliography}

\end{document}